\def\N{{\mathbb N}}
\def\C{{\mathcal C}}
\def\P{{\mathcal P}}
\documentclass{elsart}
\usepackage[cp1250]{inputenc}
\usepackage{bbm}
\usepackage{amssymb}
\usepackage{amsmath}
\usepackage{epsf}

\begin{document}

\newtheorem{theorem}{Theorem}[section]
\newtheorem{lemma}{Lemma}[section]
\newtheorem{corollary}{Corollary}[section]

\newtheorem{define}{Definition}[section]

\newcommand{\A}{\mathcal{A}}
\newcommand{\shor}{\overline{s}_{\alpha, \beta}}
\newcommand{\sdol}{\underline{s}_{\alpha, \beta}}
\newcommand{\uz}{(u_n)_{n \in \mathbb{Z}}}
\newcommand{\wn}{(w_n)_{n \in \mathbb{N}}}
\newcommand{\Q}{\mathbb{Q}[\tau]}
\newcommand{\Z}{\mathbb{Z}[\tau]}
\newcommand{\Qg}{\mathbb{Q}[\gamma]}
\newcommand{\Zg}{\mathbb{Z}[\gamma]}
\newcommand{\Qe}{\mathbb{Q}[\eta]}
\newcommand{\Ze}{\mathbb{Z}[\eta]}
\newcommand{\Qeps}{\mathbb{Q}[\varepsilon]}
\newcommand{\Zeps}{\mathbb{Z}[\varepsilon]}
\newcommand{\Zepss}{\mathbb{Z}[\varepsilon']}
\newcommand{\Qepss}{\mathbb{Q}[\varepsilon']}
\newcommand{\f}{\varphi}
\newcommand{\Si}{\Sigma}
\newcommand{\SO}{\Sigma(\Omega)}
\newcommand{\SEM}{\Sigma_{\varepsilon,\eta}(\Omega)}
\newcommand{\SE}{\Sigma_{\varepsilon,\eta}}
\newcommand{\SEE}{\Sigma_{\varepsilon,\varepsilon'}}
\newcommand{\CAP}{\Sigma_{\varepsilon,\eta}[c,c+l)}
\newcommand{\CEP}{\Sigma_{\varepsilon,\varepsilon'}[c,c+l)}
\newcommand{\CO}{\Sigma_{\varepsilon,\varepsilon'}(\Omega)}
\newcommand{\eps}{\varepsilon}
\newcommand{\sq}{\stackrel{q}{\sim}}
\newcommand{\ue}{u_{\eps,\eps'}(\Omega)}

%%%%%%%%%%%%%%%%%%%%%%%%%%%%%%%%%%%%%%%%%%%%%%%%%%%%%%%%%%%%%%%%%%%%%%
\begin{frontmatter}
\title{Factor versus palindromic complexity of uniformly recurrent infinite words}

\author{Peter Bal\'a\v zi, Zuzana Mas\'akov\'a, Edita Pelantov\'a}

\address{Department of Mathematics, FNSPE, Czech Technical University\\
Trojanova 13, 120 00 Praha 2, Czech Republic\\
e-mail: peter\_balazi@centrum.cz}

\begin{abstract}
We study the relation between the palindromic and factor
complexity of infinite words. We show that for uniformly recurrent
words one has $\P(n)+\P(n+1) \leq \Delta {\mathcal C}(n) + 2,$ for
all $n \in\N$. For a large class of words it is a better estimate
of the palindromic complexity in terms of the factor complexity
then the one presented in \cite{AllBaaCasDam}. We provide several
examples of infinite words for which our estimate reaches its
upper bound. In particular, we derive an explicit prescription for
the palindromic complexity of infinite words coding $r$-interval
exchange transformations. If the permutation $\pi$ connected with
the transformation is given by $\pi(k)=r+1-k$ for all $k$, then
there is exactly one palindrome of every even length, and exactly
$r$ palindromes of every odd length.
\end{abstract}

%\maketitle
\end{frontmatter}

%%%%%%%%%%%%%%%%%%%%%%%%%%%%%%%%%%%%%%%%%%%%%%%%%%%%%%%%%%%%%%%%%%%%%%%
\section{Introduction}

Recently, palindromes have become a popular subject of study in
the field of combinatorics on infinite words. Recall that a
palindrome is a word which remains unchanged if read backwards. In
natural language it is for example the word ``madam" in English,
or ``krk" (neck) in Czech. We shall study infinite words $u$ over
a finite alphabet $\A$, i.e. sequences $u=(u_n)_{n\in\N}$ where
$u_i\in\A$ for all $i\in\N=\{0,1,2,\dots\}$. A palindrome of the
length $n$ in the infinite word $u$ is a factor
$p=u_iu_{i+1}\cdots u_{i+n-1}$ such that $u_iu_{i+1}\cdots
u_{i+n-1}=u_{i+n-1}u_{i+n-2}\cdots u_i$.

The attractiveness of palindromes increased when Droubay and
Pirillo provided yet another equivalent definition of sturmian
words using palindromes. They have shown in~\cite{DroPir} that an
infinite word $u$ is sturmian if and only if $u$ contains exactly
one palindrome of every even length and exactly two palindromes of
every odd length.

A strong motivation for the study of palindromes in infinite words
appeared already before, in their application in modeling of solid
materials with long-range order, the so-called quasicrystals. In
1982, Dan Shechtman et al. \cite{SchBleGraCah} discovered an
aperiodic structure (which was formed by rapidly-quenched aluminum
alloys) that has icosahedral rotational symmetry, but no
three-dimensional translational invariance (see e.g.
\cite{BomTay}). The existence of such structures has been
absolutely unexpected. Since then, many other stable and unstable
aperiodic structures with crystallographically forbidden
rotational symmetry were discovered; they were named
quasicrystals.

Since the discovery of quasicrystals there has been an increasing
interest in the study of the spectral properties of non-periodic
Schr\"{o}dinger operators. One can assign to an infinite word $u$
over an alphabet $\mathcal{A}$, which models a one-dimensional
quasicrystal, a Schr\"{o}dinger operator $H$ acting on the Hilbert
space $\ell^2(\mathbb{Z})$ as follows
$$
(H \phi)(n) = \phi(n+1) + \phi(n-1) + V(u_n) \phi (n),
$$
where $V : \mathcal{A} \mapsto \mathbb{R}$ is an injection and
represents a potential of the operator.

Many nice properties of these operators have been shown, and they
are well understood at least in the one-dimensional case. The
survey papers \cite{Da}, \cite{Su} map the history of this effort.
One of the main tasks is to derive the spectral properties of the
Schr\"{o}dinger operator $H$ from the properties of the sequence
$V(u_n)$.
The physical motivation behind this study is that the spectral
properties of operators determine the conductivity of the given
structure. Very roughly speaking, if the spectrum is pure point
then the structure is behaving like an insulant. In case of
absolutely continuous spectrum the material is becoming a
conductor.

%In the beginning, the studies of aperiodic Schr\"{o}dinger
%operators used a potential, given by the simplest aperiodic
%structure, the Fibonacci word. It was shown that the spectrum of
%such an operator is purely singular continuous with zero measure.
%This fact reflects that the operator models structures which are
%in between periodic (corresponding spectrum is absolutely
%continuous) and disordered (pure point spectrum).

Generally, the task of describing the spectral properties of the
operator with potential given by an arbitrary infinite word $u$ is
not a simple one. The relevance of the study of palindromes in the
infinite words has been proven by Hof et al. \cite{HoKnSi} who
showed that the operators given by words having arbitrary large
palindromes have purely singular continuous spectrum.

The aim of this article is to find a relation between factor and
palindromic complexity of uniformly recurrent words. Let us first
introduce the basic notions which will be used in sequel.

The set of all factors of length $n$ of an infinite word
$u=u_0u_1u_2\cdots$ is denoted by
$$
{\mathcal L}_n(u) = \{ w_1\cdots w_n \mid \exists i\in\N,\
w_1\cdots w_n = u_i\cdots u_{i+n-1}\}\,.
$$
The set of all factors of $u$, including the empty word
$\varepsilon$ is called the language of $u$ and denoted
$$
{\mathcal L}(u)= \bigcup_{n\in\N}{\mathcal L}_n(u)\,.
$$
The variability of local configurations in the word $u$ is
characterized by the factor complexity, the function $\C:\N\to\N$,
given by the prescription
$$
\C(n):=\# {\mathcal L}_n(u)\,.
$$
It is known that if there is an $n_0$ such that $\C(n_0)\leq n_0$,
then the word $u$ is eventually periodic, i.e.\ there exists
$k\in\N$ such that $u_{k+n}=u_n$ for every $n\geq n_0$. Any
aperiodic (i.e. not eventually periodic) word therefore satisfies
$\C(n)\geq n+1$ for every $n\in\N$. Aperiodic words of minimal
complexity $\C(n)=n+1$ are called sturmian words. For a survey of
different characteristics and properties of sturmian words
see~\cite{berstel}.

The mirror image, or reversal, of a finite word $w=w_1\cdots w_n$
is the word $\overline{w}=w_n\cdots w_1$. If the language
${\mathcal L}(u)$ contains with every factor $w$ also its mirror
image $\overline{w}$, we say that ${\mathcal L}(u)$ is invariant
under reversal.

The palindromic complexity of the infinite word
$u=(u_{n})_{n\in\N}$ is a function $\P: \N\to\N$ which counts the
number of palindromes of a given length. Formally,
$$
\P(n) := \# \{ w\in {\mathcal L}_n(u) \mid w=\overline{w}\}\,.
$$

Trivially, one has $\P(n)\leq \C(n)$. A non-trivial result is an
estimate of $\P(n)$ using $\C(n)$ provided in~\cite{AllBaaCasDam}.

%\bigskip
\begin{theorem}[\cite{AllBaaCasDam}]\label{thm:ABCD}
For arbitrary infinite word one has
\begin{equation}\label{eq:abcd}
\P(n)\leq \frac{16}{n}\, \C\bigl( n + \lfloor\frac{n}{4}
\rfloor\bigr)\,,\qquad \hbox{ for all }\ n\in\N\,.
\end{equation}
\end{theorem}

Let us mention that Theorem~\ref{thm:ABCD} implies the result
of~\cite{DaZare}: The palindromic complexity of a fixed point of a
primitive morphism is bounded.

In this paper we provide an estimate of $\P(n)$ of uniformly
recurrent words using the first difference $\Delta\C(n) :=
\C(n+1)-\C(n)$. For words whose factor complexity is a polynomial
of degree $\leq 16$, this estimate is better than that
of~\eqref{eq:abcd}. Let us recall that an infinite word $u$ is
uniformly recurrent, if the gaps between consecutive occurrences
of any factor $w\in{\mathcal L}(u)$ in the word $u$ are bounded.
Equivalently, if for every $n\in\N$ there exists $R(n)\in\N$ such
that in an arbitrary segment of length $R(n)$ in the word $u$ one
finds all factors of ${\mathcal L}_n(u)$, i.e.
$$
{\mathcal L}_n(u) = \{u_iu_{i+1}\cdots u_{i+n-1} \mid k\leq i\leq
k+R(n) \}\,,\qquad \hbox{ for all }\ k,n\in\N\,.
$$

Let us mention that sturmian words are example of uniformly
recurrent words with language closed under reversal~\cite{DroPir}.

In section~\ref{sec:odhad} we show the following theorem.

%\bigskip
\begin{theorem}\label{thm:odhad}
Let $u=(u_n)_{n\in\N}$ be an uniformly recurrent word.
\begin{itemize}
\item[(i)] If ${\mathcal L}(u)$ is not closed under reversal, then
$\P(n)=0$ for sufficiently large $n$.
\item[(ii)] If ${\mathcal L}(u)$ is closed under reversal, then
$$
\P(n)+\P(n+1) \leq \Delta {\mathcal C}(n) + 2 \,,\qquad\hbox{ for
all }\ n\in\N\,.
$$
\end{itemize}
\end{theorem}

It is interesting that equality in the latter estimate of the
palindromic complexity holds for some known classes of infinite
words, such as Arnoux-Rauzy words or fixed points of canonical
substitutions associated to numeration systems with base $\beta$,
where $\beta$ is a Parry number~\cite{FrMaPe}. We list some of
these examples in section~\ref{sec:odhad}. In
section~\ref{sec:rIET} we show that the equality in the estimate
is valid also for infinite words coding $r$-interval exchange
transformation.

%%%%%%%%%%%%%%%%%%%%%%%%%%%%%%%%%%%%%%%%%%%%%%%%%%%%%%%%%%%%%%%%%%%%%%%
\section{Proof of Theorem~\ref{thm:odhad}}\label{sec:odhad}

First we show that unboundedness of the length of palindromes in
an infinite uniformly recurrent word $u$ implies that the language
of $u$ is invariant under mirror image.

\begin{lemma}
Let $u$ be an infinite word which is uniformly recurrent and such
that $\limsup_{n\to\infty}\P(n)>0$. Then $\overline{{\mathcal
L}(u)} = {\mathcal L}(u)$.
\end{lemma}

\pf Let $n\in\N$. Consider $R(n)$ from the definition of uniformly
recurrent words. Let $p$ be a palindrome of length greater than
$R(n)$. It contains all factors of $u$ of length $n$. In the same
time it contains with every factor $w$ also its mirror image. Thus
$\overline{{\mathcal L}_n(u)} = {\mathcal L}_n(u)$ for all
$n\in\N$.

The above lemma in fact proves (i) of Theorem~\ref{thm:odhad}.
Crucial tool for the proof of (ii) is the notion of a Rauzy graph
of an infinite word.

Let $u=(u_n)_{n\in\N}$ be an infinite word, $n\in\N$. The Rauzy
graph $\Gamma_n$ of $u$ is an oriented graph whose set of vertices
is ${\mathcal L}_n(u)$ and the set of edges is ${\mathcal
L}_{n+1}(u)$. An edge $e\in{\mathcal L}_{n+1}(u)$ starts at the
vertex $x$ and ends at the vertex $y$, if $x$ is a prefix and $y$
is a suffix of $e$.

\begin{figure}[ht]
\begin{center}
\begin{picture}(220,30)
\put(50,20){\circle*{5}} \put(210,20){\circle*{5}}
\put(60,20){\vector(1,0){140}}
%\qbezier(50,20)(120,40)(190,20)
\put(10,8){$x=w_0w_1\cdots w_{n-1}$} \put(170,8){$y=w_1\cdots
w_{n-1}w_n$} \put(82,25){$e=w_0w_1\cdots w_{n-1}w_n$}
\end{picture}
\end{center}
\caption{Incidence relation between an edge and vertices in a
Rauzy graph.} \label{f}
\end{figure}
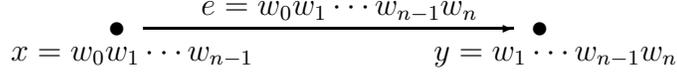

If the word $u$ is uniformly recurrent, the graph $\Gamma_n$ is
strongly connected for every $n\in\N$, i.e. there exists an
oriented path from every vertex $x$ to every vertex $y$ of the
graph.

The outdegree of a vertex $x\in{\mathcal L}_n(u)$ is the number of
edges which start in $x$. It is denoted by ${\rm deg}_+(x)$,
$$
{\rm deg}_+(x) := \#\{ a\in\A \mid xa\in {\mathcal
L}_{n+1}(u)\}\,.
$$
Similarly, we define the indegree of $x$ as
$$
{\rm deg}_-(x) := \#\{ a\in\A \mid ax\in {\mathcal
L}_{n+1}(u)\}\,.
$$

The sum of outdegrees over all vertices is equal to the number of
edges in every oriented graph. Similarly, it holds for indegree.
In particular, for the Rauzy graph we have
$$
\sum_{x\in{\mathcal L}_{n}(u)}{\rm deg}_+(x) = \# {\mathcal
L}_{n+1}(u) = \sum_{x\in{\mathcal L}_{n}(u)}{\rm deg}_-(x)\,.
$$
Since $\Delta\C(n) = \# {\mathcal L}_{n+1}(u) - \# {\mathcal
L}_{n}(u)$, we obtain
\begin{equation}\label{eq:delta}
\Delta\C(n) = \sum_{x\in{\mathcal L}_{n}(u)}\bigl({\rm deg}_+(x) -
1 \bigr) = \sum_{x\in{\mathcal L}_{n}(u)}\bigl({\rm deg}_-(x) - 1
\bigr)\,.
\end{equation}
A non-zero contribution to $\Delta\C(n)$ is therefore given only
by those factors $x\in{\mathcal L}_n(u)$, for which ${\rm
deg}_+(x)\geq 2$, i.e. such that there exist distinct letters
$a,b\in\A$ satisfying $xa,xb\in{\mathcal L}_{n+1}(u)$. A factor of
$u$, which has at least two extensions to the right is called a
right special factor of $u$. Similarly one can define a left
special factor, and the relation~\eqref{eq:delta} can be rewritten
as
$$
\Delta\C(n) = \sum_{x\in{\mathcal L}_{n}(u),\ \hbox{\scriptsize
$x$ right special}\hspace*{-1.5cm}}\hspace*{-0.3cm}\bigl({\rm
deg}_+(x) - 1 \bigr) = \sum_{x\in{\mathcal L}_{n}(u),\
\hbox{\scriptsize $x$ left
special}\hspace*{-1.5cm}}\hspace*{-0.3cm}\bigl({\rm deg}_-(x) - 1
\bigr)\,.
$$

\vskip0.4cm \noindent {\bf PROOF of (ii) of
Theorem~\ref{thm:odhad}.} Suppose  that the language of the
infinite word $u$ is closed under reversal. Consider the operation
$r$ which to every vertex of the graph associates
$\rho(x)=\overline{x}$ and to every edge associates
$\rho(e)=\overline{e}$.

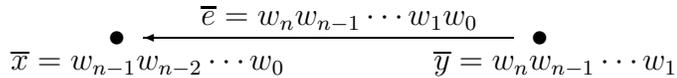
\begin{figure}[ht]
\begin{center}
\begin{picture}(220,30)
\put(50,20){\circle*{5}} \put(210,20){\circle*{5}}
\put(200,20){\vector(-1,0){140}}
%\qbezier(50,20)(120,40)(190,20)
\put(10,8){$\overline{x}=w_{n-1}w_{n-2}\cdots w_{0}$}
\put(170,8){$\overline{y}=w_{n}w_{n-1}\cdots w_1$}
\put(82,25){$\overline{e}=w_nw_{n-1}\cdots w_{1}w_0$}
\end{picture}
\end{center}
\caption{The action of the mapping $\rho$ on the edge and the
vertices of Figure~\ref{f}.} \label{f2}
\end{figure}

This operation maps the Rauzy graph $\Gamma_n$ onto itself.
Obviously,
$$
\begin{array}{rcl}
\P(n) &=& \#\{ z\in{\mathcal L}_n(u) \mid \rho(z)=z\}\,,\\
\P(n+1) &=& \#\{ e\in{\mathcal L}_{n+1}(u) \mid \rho(e)=e\}\,.
\end{array}
$$
We shall be interested in the pathes leading between special
factors. More precisely, we shall call a {\it simple path} an
oriented path $w=v_0v_1\ldots v_k$, such that its initial vertex
$v_0$, and its final vertex $v_k$ are left or right special
factors, and the other vertices are not special factors, i.e.
${\rm deg}_+(v_i)={\rm deg}_-(v_i)=1$ for $i=1,2,\dots,k-1$. A
special factor is considered as a simple path of length~$0$. Since
the infinite word $u$ is uniformly recurrent, the graph $\Gamma_n$
is strongly connected, and therefore every vertex and every edge
belongs to a simple path.

For an edge $e$ satisfying $\rho(e)=e$ we find the simple path $w$
which contains $e$. Since $\rho(e)=e$, the operation $r$ must map
the path $w$ onto itself. Similarly, if for a vertex $z$ it holds
that $\rho(z)=z$, then the  simple path containing $z$ is mapped
by $\rho$ onto itself.

To give an upper bound on $\P(n)+\P(n+1)$ therefore consists in
finding the number of simple paths in the Rauzy graph $\Gamma_n$
which are mapped by $\rho$ onto itself. It therefore suffices to
study the so-called reduced Rauzy graph.

The set $V$ of vertices of the reduced Rauzy graph is formed by
all $x\in{\mathcal L}_n(u)$ which are either left or right special
factors of $u$. Two vertices $x,y\in V$ are connected by an
oriented edge from $x$ to $y$, if in the original Rauzy graph
$\Gamma_n$ there exists a simple path from $x$ to $y$. The
operation $\rho$ maps the reduced Rauzy graph onto itself.

The set $V$ of vertices of the reduced Rauzy graph can be divided
into disjoint cycles of the mapping $\rho$. Since $\rho^2={\rm
Id}$, the cycles are either of length $1$ or $2$. The cycles of
length $1$ are given by special factors invariant under $\rho$,
i.e.\ special factors, which are themselves palindromes. Let us
denote their number by $\alpha$, and denote the number of cycles
of length $2$ by $\beta$. Note that the number of vertices in the
reduced Rauzy graph (i.e.\ left or right special factors in
$\Gamma_n$) is $\alpha+2\beta$.

If there is an edge from $x$ to $y$, where $x$ and $y$ belong to
different cycles, then there is another edge leading from
$\rho(y)$ to $\rho(x)$. Since the reduced Rauzy graph is strongly
connected, the number of edges, which lead between vertices of
different cycles, is at least $2(\alpha+\beta-1)$. These edges
correspond in the original Rauzy graph $\Gamma_n$ to the simple
paths of non-zero length which are not mapped by $\rho$ onto
itself.

As we have said, the number of palindromes of length $n$ and $n+1$
is bounded by the number of simple paths in $\Gamma_n$, which are
mapped by $\rho$ onto itself. We thus have
$$
\P(n) + \P(n+1) \leq \sum_{\hbox{\scriptsize $x$ is left or right
special}\hspace*{-1.5cm}}\hspace*{-0.3cm}{\rm deg}_+(x) \ -\
2(\alpha+\beta-1) \ +\ \alpha\,,
$$
where the first summand is the number of all simple paths of
non-zero length in $\Gamma_n$, the second summand estimates the
number of simple paths of non-zero length which are not mapped
onto itself, and the third one is the number of simple paths of
zero length invariant under $\rho$ (i.e.\ palindromic special
factors). We obtain
$$
\begin{aligned}
\P(n) + \P(n+1)
 &\leq \sum_{\hbox{\scriptsize $x$ is left or right
special}\hspace*{-1.5cm}}\hspace*{-0.3cm}{\rm deg}_+(x) \ -\
(\alpha+2\beta) \ +\
2 \ = \\[4mm]
&=\sum_{\hbox{\scriptsize $x$ is left or right
special}\hspace*{-1.5cm}}\hspace*{-0.3cm}\bigl({\rm deg}_+(x) - 1
\bigr) \ + \ 2 \ = \ \Delta\C(n) \ +\ 2\,,
\end{aligned}
$$
where we have used that $\alpha+2\beta$ is the number of left or
right special factors in $\Gamma_n$. This completes the proof.

%%%%%%%%%%%%%%%%%%%%%%%%%%%%%%%%%%%%%%%%%%%%%%%%%%%%%%%%%%%%%%%%%%%%%%%
\section{Examples of infinite words with maximal $\P(n)+\P(n+1)$}\label{sec:examples}

In this section we present several examples of infinite words
which satisfy $\limsup_{n\to\infty}\P(n)>0$ and
\begin{equation}\label{eq:rovnost}
\P(n)+\P(n+1)= \Delta\C(n)+2\,.
\end{equation}
These are in a sense words with maximal number of palindromes.

\paragraph*{1. Arnoux-Rauzy sequences.} Arnoux-Rauzy sequences are generalizations of sturmian
words for an alphabet with more than 2 letters. An infinite word
$u$ over an $r$-letter alphabet is called Arnoux-Rauzy of order
$r$, if for every $n\in\N$ there exists exactly one left special
factor, say $w_L$, and exactly one right special factor, say
$w_R$, of length $n$, and they satisfy ${\rm deg}_+(w_R)={\rm
deg}_-(w_L)=r$. Note that Arnoux-Rauzy sequences of order $2$ are
precisely the sturmian words. Directly from the definition one can
deduce that the factor complexity is $\C(n)=(r-1)n+1$ for all
$n\in\N$. In~\cite{DaZa} it was shown that
$$
\P(n)=\biggl\{\begin{array}{cl}
 r,&\quad \hbox{if $n$ is odd,}\\
 1,&\quad \hbox{if $n$ is even.}
 \end{array}%\right.
$$
Since $\Delta\C(n) = r-1$, we obtain
$$
\P(n)+\P(n+1) = r+1 = \Delta\C(n)+2\,,
$$
and thus the Arnoux-Rauzy words satisfy~\eqref{eq:rovnost}.

\paragraph*{2. Complementation-symmetric
sequences.} In~\cite{AllBaaCasDam} it was shown that
complementation-symmetric sequences, with factor complexity
$\C(n)=2n$ for all $n\in\N$, satisfy $\P(n)=2$ for all $n\geq 1$.
Recall that a complementation-symmetric sequence on a two-letter
alphabet, say $\A=\{a,b\}$, is a sequence such that for any factor
occurring in it, the word obtained by changing $a$'s into $b$'s
and vice versa, is also a factor. Since $\Delta\C(n)=2$, we have
again
$$
\P(n)+\P(n+1) = 4 = \Delta\C(n)+2\,.
$$

\paragraph*{3. Words associated with $\beta$-integers.}
In~\cite{palindromy} one studies palindromes in words associated
with $\beta$-integers, i.e. positive real numbers which have
vanishing fractional part in the numeration system with base
$\beta$. For the description of the words $u_\beta$ we introduce
the R\'enyi expansion of $1$.

Let $\beta$ be a fixed real number, $\beta>1$. Denote by $T_\beta$
the mapping $T_\beta:[0,1]\to[0,1)$, given by the prescription
$$
T_\beta(x):=\beta x-\lfloor \beta x \rfloor\,.
$$
The sequence
$$
d_\beta(1)=t_1t_2t_3\cdots\,,\quad\hbox{ where }\quad t_i:=\lfloor
T_\beta^{i-1}(1) \rfloor\,,\quad i=1,2,3,\dots
$$
is called the R\'enyi expansion of $1$. If $d_\beta(1)$ is
eventually periodic, then $\beta$ is called a Parry number.

The infinite word $u_\beta$, which codes the distances between
$\beta$-integers, is the fixed point of a morphism over a finite
alphabet. The morphisms are of two types, according to the type of
the Parry number $\beta$.

\begin{itemize}
\item If $d_\beta(1)=t_1\cdots t_m0^\omega$, with $t_m\neq 0$,
then $\beta$ is called a simple Parry number. In this case
$u_\beta$ is the fixed point of the substitution
$\varphi=\varphi_\beta$ over the alphabet $\A=\{0,1,\cdots,m-1\}$,
given by
$$
\begin{array}{rcl}
\varphi(0)&=&0^{t_1}\ 1,\\
\varphi(1)&=&0^{t_2}\ 2,\\
&\vdots&\\
\varphi({m\!-\!2})&=&0^{t_{m-1}}\ (m\!-\!1),\\
\varphi({m\!-\!1})&=&0^{t_m}.
\end{array}
$$

\item
If $d_\beta(1)=t_1\cdots t_m(t_{m+1}\cdots t_{m+p})^\omega$, where
$m$, $p$ are minimal indices which allow such notation, then
$u_\beta$ is the fixed point of the substitution
$\varphi=\varphi_\beta$ over the alphabet
$\A=\{0,1,\cdots,m+p-1\}$, given by
$$
\begin{array}{rcl}
\varphi(0)&=&0^{t_1}\ 1,\\
\varphi(1)&=&0^{t_2}\ 2,\\
&\vdots&\\
\varphi(m\!-\!1)&=&0^{t_m}\ m,\\
&\vdots&\\
\varphi({m\!+\!p\!-\!2})&=&0^{t_{m+p-1}}\ (m\!+\!p\!-\!1),\\
\varphi({m\!+\!p\!-\!1})&=&0^{t_{m+p}}\ m.
\end{array}
$$
\end{itemize}

For infinite word $u_\beta$ one can easily show that are uniformly
recurrent. The condition of invariance of the language of
$u_\beta$ under reversal is described in~\cite{FrMaPe} for the
case of simple Parry number, i.e. $d_\beta(1)=t_1\cdots
t_m0^\omega$. It is shown that ${\mathcal L}(u_\beta)$ is closed
under reversal if and only if $t_1=\cdots =t_{m-1}\geq t_m$. For
the case $d_\beta(1)=t_1\cdots t_m(t_{m+1}\cdots t_{m+p})^\omega$
it is shown in~\cite{bernat} that the language of $u_\beta$ is
closed under reversal if and only if $m=p=1$. Papers~\cite{FrMaPe}
and~\cite{Lubka} show that if $u_\beta$ has the language invariant
under reversal, then
$$
\P(n+2)-\P(n) = \Delta^2\C(n)=\C(n+1) - \C(n)\,,
$$
which allows one to derive the validity of~\eqref{eq:rovnost}.

While in examples 1 and 2 the second difference
$\Delta^2\C(n)\equiv 0$, for the words $u_\beta$ it holds that
$\Delta^2\C(n)\in\{-1,0,1\}$ and all three values are reached
infinitely many times.

\paragraph*{4. Words coding $r$-interval exchange transformation.}
Another possible generalization of sturmian words are words coding
a bijective transformation of the interval $[0,1)$ onto itself,
known under the name $r$-interval exchange. Let us recall the
definition of an interval exchange map. It can be found together
with some properties in \cite{DroPir}, \cite{Keane}.

Given $r$ positive numbers $\alpha_1, \alpha_2, \ldots, \alpha_r$
such that $\sum^{r}_{i =1} \alpha_i = 1$. They define a partition
of the interval $I = [0,1)$ into $r$ intervals
$$
I_k = \biggl[\,\sum^{k-1}_{i =1} \alpha_i, \ \sum^{k}_{i =1}
\alpha_i\biggr)\,,\qquad k = 1,2, \ldots, r\,.
$$
Let $\pi$ denote a permutation of the set $\{1,2,\ldots, r\}$. The
interval exchange transformation associated with $\alpha_1,\ldots,
\alpha_r$ and $\pi$ is defined as the map $T:I\to I$ which
exchanges the intervals $I_k$ according to the permutation $\pi$,
$$
T(x) = x + \sum_{j < \pi(k)} \alpha_{\pi^{-1}(j)} - \sum_{j < k}
\alpha_j\,, \quad \textrm{ for } x \in I_k\,.
$$
For $x_0 \in I$, the sequence $(T^n (x_0))_{n \in \mathbb{Z}}$ is
called the orbit of $x_0$ under $T$. The infinite bidirectional
word $\uz$ over the alphabet $\mathcal{A} = \{1, \ldots, r\}$
associated to the orbit $(T^n(x_0))_{n \in \mathbb{Z}}$ is defined
as
$$
u_n = k \in \mathcal{A} \quad \Leftrightarrow \quad  T^n(x_0) \in
I_k\,.
$$
The complexity of the word corresponding to any $r$-interval
exchange transformation satisfies $\mathcal{C}(n) \leq n(r-1) +
1$, for all $n\in{\mathbb N}$. Here we focus on the
non-degenerated case, i.e.\ on mappings $T$ for which the
complexity of the word associated to the orbit of arbitrary $x_0
\in I$ satisfies $\mathcal{C}(n) = (r-1)n + 1$, for all $n \in
\mathbb{N}$. This property is ensured by additional conditions
(denoted by $\mathfrak{P}$) on the parameters of the map $T$.
$$
(\mathfrak{P})\qquad
\begin{array}{cl}
1.& \alpha_1, \ldots, \alpha_r \hbox{ are linearly independent over } \mathbb{Q},\\[2mm]
2.& \pi \{ 1, \ldots, k \} \neq \{ 1, \ldots, k \} \hbox{ for each
} k = 1,2, \ldots, r-1.
\end{array}
$$
If the conditions ($\mathfrak{P}$) are fulfilled, then the set
$\{T^n(x_0)\}_{n \in \mathbb{Z}}$ is dense in $I$ for each $x_0
\in I$ and the dynamical system associated to the transformation
$T$ is minimal. It implies that the infinite word corresponding to
the sequence $(T^n(x_0))_{n \in \mathbb{Z}}$ is uniformly
recurrent.

Another important consequence of ($\mathfrak{P}$) is that the
language of the word $\uz$ corresponding to $(T^n(x_0))_{n \in
\mathbb{Z}}$ does not depend on the position of the starting point
$x_0$, but only on the transformation $T$. Therefore the notation
$\mathcal{L}(T)$, which we adopt here, is justified. We know that
${\mathcal L}(T_1)={\mathcal L}(T_2)$ only if $T_1$ and $T_2$
coincide.

If $r=2$, the permutation satisfying $\mathfrak{P}$ is $\pi(1)=2$,
$\pi(2)=1$ and the corresponding word is sturmian. On the other
hand, every sturmian word can be obtain as a coding of a
2-interval exchange transformation.

If $r=3$, then the condition 2. of $\mathfrak{P}$ is satisfied by
three permutations. One can easily see that only the permutation
$\pi(1)=3$, $\pi(2)=2$, $\pi(3)=1$ gives an infinite word with
language invariant under reversal. Such words can be geometrically
represented by cut-and-project sequences~\cite{GuMaPe}.

For general $r$, the language of the infinite word $u$ closed
under reversal if and only if
\begin{equation}\label{cislo3}
\pi(1) = r, \quad \pi(2) = r -1, \quad \ldots, \quad \pi(r) = 1.
\end{equation}
Only for such permutation the infinite word $u$ coding the
corresponding interval exchange transformation one may have
$\limsup_{n\to\infty}\P(n)>0$.

The palindromic complexity in words coding 3-interval exchange map
was described in~\cite{DaZa}. In section~\ref{sec:rIET} we
generalize their result for any $r$. We show that for words coding
an $r$-interval exchange transformation with
permutation~\eqref{cislo3} the equality~\eqref{eq:rovnost} holds.

%%%%%%%%%%%%%%%%%%%%%%%%%%%%%%%%%%%%%%%%%%%%%%%%%%%%%%%%%%%%%%%%%%%%%%%
\section{Words Coding Interval Exchange Transformation}\label{sec:rIET}

In this section we will be dealing only with such transformations
$T$ of $r$-intervals for which the permutation $\pi$ satisfies
(\ref{cislo3}). In this case the transformation has the form of
\begin{equation}\label{eq:tri}
T(x) = x + \sum_{j > k} \alpha_{j} - \sum_{j < k} \alpha_j \ \
\textrm{ for } x \in I_k.
\end{equation}
It is known that there exists an interval $I_w \subset I_{w_{0}}$
for every word $w = w_0 w_1 \ldots w_{n-1} \in \mathcal{L}(T)$
such that the sequence of points $x, T(x), \ldots, T^{n-1}(x)$ is
coded by the same word $w$ for each $x \in I_w$. Note that the
boundaries of the interval $I_w$ belong to the set
$\mathbb{Z}[\alpha_1, \ldots, \alpha_r] = \{ \sum k_i \alpha_i \ |
\ k_i \in \mathbb{Z} \}$.

Let us denote the decomposition of the interval $I = [0,1)$ by the
transformation $T^{-1}$ by $\tilde{I}_1, \tilde{I}_2, \ldots,
\tilde{I}_r$ and analogously $\tilde{I}_w$ for an arbitrary $w \in
\mathcal{L}(T^{-1})$.

Clearly, $\tilde{I}_{\pi^{-1}(j)} = T(I_j)$ for each $j \in \{1,
\ldots, r\}$. Since $\pi$ is of the form~\eqref{cislo3}, it
follows that  $I_j=[a,b)$ implies $\tilde{I}_{\pi^{-1}(j)} =
T(I_j)=[1-b,1-a)$. The same relation is therefore valid for any
factor $w\in{\mathcal L}(T)$,
\begin{equation}\label{eq:nevim}
I_w=[a,b) \quad\implies\quad \tilde{I}_{\pi^{-1}(w)} =[1-b,1-a)\,.
\end{equation}

Now we have everything prepared for determination of the
palindromic complexity.

\begin{theorem}
Let $\alpha_1, \ldots, \alpha_r$ be positive real numbers,
linearly independent over $\mathbb{Q}$ and $\pi$ a permutation
satisfying (\ref{cislo3}). Then
$$
\mathcal{P}(n) = \left\{
\begin{array}{ll} 1 & \textrm{ for each } n \textrm{ even }, \\
                 r & \textrm{ for each } n \textrm{ odd }.
\end{array} \right.
$$
\end{theorem}

%\begin{proof}
\pf Consider the palindrome of even length in the form of
$$
w_{n-1} w_{n-2} \ldots w_0 w_0 \ldots w_{n-2} w_{n-1} \in
\mathcal{L}(T).
$$
It means that there exists $x \in [0,1)$ such that
$$
\begin{array}{rlrlllrl}
 x \in I_{w_0}\,, & &  T(x) \in I_{w_1}\,, &  & \ldots\,, &  &  T^{n-1}(x) \in I_{w_{n-1}}\,, \\
 T^{-1}(x) \in I_{w_0}\, ,& & T^{-2}(x) \in I_{w_1}\, ,&  & \ldots\,, &
& T^{-n}(x) \in I_{w_{n-1}}\,.
\end{array}
$$
Hence  $x \in I_w$, where $w = w_0, \ldots w_{n-1}$ and on the
other side
$$
\begin{array}{rclllll} x & \ \in\ & T(I_{w_0}) =  \tilde{I}_{\pi^{-1}(w_0)}, \\
                      T^{-1}(x) & \ \in\  & T(I_{w_1}) =  \tilde{I}_{\pi^{-1}(w_1)}, \\
                      & \vdots & \\
                      T^{-n+1}(x) & \ \in\  & T(I_{w_{n-1}}) = \tilde{I}_{\pi^{-1}(w_{n-1})}.
\end{array}
$$
It follows that $x \in \tilde{I}_{\pi^{-1}(w_0 w_1 \ldots
w_{n-1})}$. Thus $x$ has to belong to the intersection of both
intervals, i.e. $x \in I_w \cap \tilde{I}_{\pi^{-1}(w)}$. If
$I_w=[a,b)$, then according to~\eqref{eq:nevim}
$$
x\in[a,b)\cap[1-b,1-a)\,.
$$

Now we use a simple fact that for every interval $[a,b)$ it holds
that
\begin{equation}\label{vim}
[a,b)\cap[s-b,s-a)\neq  \emptyset \quad \Longleftrightarrow \quad
\frac{s}{2} \in [a,b)\,.
\end{equation}
Therefore
\begin{equation}
\frac{1}{2} \in I_w \cap \tilde{I}_{\pi^{-1}(w)}.
\end{equation}
We have shown that every palindrome of even length arises from the
coding of
$$
T^{-n} \left( \tfrac{1}{2} \right), \ldots, T^{-1} \left(
\tfrac{1}{2} \right), \tfrac{1}{2}, T \left( \tfrac{1}{2} \right),
\ldots, T^{n -1} \left( \tfrac{1}{2} \right).
$$

The fact that $\{T^n (x)\}_{n \in \mathbb{Z}}$ is dense in $[0,1)$
implies that the previous sequence occurs in $(T^n(x))_{n \in
\mathbb{Z}}$, for each  $x$. Thus the coding of $(T^n(x))_{n \in
\mathbb{Z}}$ includes exactly one palindrome of even length for
each $n$.

Consider now the palindrome of odd length in the form of
$$
w_{n-1} w_{n-2} \ldots w_1 w_0 w_1 \ldots w_{n-2} w_{n-1} \in
\mathcal{L}(T).
$$
Again, it means that there exist $x,y \in [0,1)$ such that
$$
\begin{array}{rrlllll} x \in I_{w_0}\,, &   T(x) \in I_{w_1}\,,  & \ \ldots \,,& \  T^{n-1}(x) \in I_{w_{n-1}}\,, \\
                y \in I_{w_1}\,,   & T^{-1}(y) \in I_{w_0} \,, &  \ \ldots \,, & \ T^{-n}(y) \in I_{w_{n-1}}\,.
\end{array}
$$
The first sequence is the coding of the word $w = w_0 w_1 \ldots
w_{n-1}$, i.e. $x \in I_w$, and the following is true for the
second one
$$
\begin{array}{rclllll} y & \ \in\ & T(I_{w_0}) =  \tilde{I}_{\pi^{-1}(w_0)}, \\
                       T^{-1}(y) & \ \in\  & T(I_{w_1}) =  \tilde{I}_{\pi^{-1}(w_1)}, \\
                       & \vdots & \\
                       T^{-n+1}(y) & \ \in\ & T(I_{w_{n-1}}) = \tilde{I}_{\pi^{-1}(w_{n-1})}.
\end{array}
$$
Thus $y \in \tilde{I}_{\pi^{-1}(w)}$. If there exists a palindrome
of odd length with the central letter $w_0$ then it has to be $y =
T(x) = x + s_{w_0}$, where $s_{w_0} \neq 0$ is a shift of $x \in
I_{w_0}$ by the mapping $T$. Using \eqref{eq:tri} we have $s_{w_0}
= \sum_{j>w_0}\alpha_j - \sum_{j<w_0}\alpha_j$. In other words we
have $x \in I_w$ and $x \in \tilde{I}_{\pi^{-1}(w)} - s_{w_0}$. If
$I_w=[c,d)$ then $\tilde{I}_{\pi^{-1}(w)}= [1-d,1-c)$ and
therefore $x\in [c,d) \cap [1-d -s_{w_0},1-c -s_{w_0})$. According
to \eqref{vim}
$$
x_{w_0}:=\frac{1 - s_{w_0}}{2} \in I_w \cap
(\tilde{I}_{\pi^{-1}(w)} - s_{w_0}).
$$
We have shown that the palindrome $w_{n-1} w_{n-2} \ldots w_1 w_0
w_1 \ldots w_{n-2} w_{n-1}$ can be obtained by
 coding of following sequences
\begin{equation}\label{sequence}
T^{-n+1} \left( x_{w_0} \right), \ldots, T^{-1} \left( x_{w_0}
\right), x_{w_0}, T \left( x_{w_0} \right), \ldots, T^{n -1}
\left( x_{w_0} \right).
\end{equation}
One may rewrite
$$x_{w_0} = \frac{1 - s_{w_0}}{2} = \frac{\sum_{j=1}^{r} \alpha_j -
\sum_{j=w_0+1}^{r} \alpha_j + \sum_{j=1}^{w_0-1} \alpha_j}{2} =
\sum_{j=1}^{w_0-1} \alpha_j + \frac{\alpha_{w_0}}{2}.
$$
It means that the point $x_{w_0}$, which correspond to central
letter $w_0$ in the palindrome of odd length, is laying in the
middle of interval $I_{w_0}$ associated to the central letter.

On the other hand, if $x_{w_0}$ is the center of one of the
intervals $I_1, \ldots, I_r$, the sequence \eqref{sequence}
corresponds to a palindrome.  Therefore  ${\P}(2n+1) = r$.

%\end{remark}

Note that according to the previous theorem the interval exchange
transformation with a permutation satisfying (\ref{cislo3}) has
the same palindromic complexity and also factor complexity as
Arnoux-Rauzy words over $r$ letters \cite{DaZa}, \cite{JusPir}.

%%%%%%%%%%%%%%%%%%%%%%%%%%%%%%%%%%%%%%%%%%%%%%%%%%%%%%%%%%%%%%%%%%%%%%%
\section{Conclusions}\label{sec:zaver}

The main result of this paper is the estimate of the palindromic
complexity of infinite words in terms of their factor complexity.
We have shown in Theorem~\ref{thm:odhad} that uniformly recurrent
words with infinitely many palindromes satisfy the following
relation
$$
\P(n)+\P(n+1) \leq \Delta {\mathcal C}(n) + 2 \,,\qquad\hbox{ for all }\ n\in\N\,.
$$

It is interesting to mention that the first difference of factor
complexity was already useful for estimation of the frequencies of
factors. In~\cite{Bosh} it is shown that the frequencies of
factors of length $n$ in a recurrent word take at most $3\Delta
\C(n)$ values.

The second part of the paper is devoted to infinite words for
which $\P(n)+\P(n+1)$ in Theorem~\ref{thm:odhad} reaches the upper
bound. We cite several examples of such infinite words among the
words for which the palindromic and factor complexity was known.
As a new result, we derive the palindromic complexity for infinite
words coding $r$-interval exchange transformation and prove that
for this class of infinite words the equality in the estimate
hold, too.

According to our knowledge all known examples of infinite words
which satisfy the equality $\P(n)+\P(n+1) = \Delta {\mathcal C}(n)
+ 2 $ for $n\in\N$ have sublinear factor complexity. A known
example of an infinite word with higher factor complexity are the
billiard sequences on three letters, for which $\C(n) = n^2+n+1$.
As shown in \cite{Borel}, they satisfy $\P(n) + \P(n+1) = 4$, and
thus billiard sequences do not reach the upper bound in
Theorem~\ref{thm:odhad}.

The proof of Theorem~\ref{thm:odhad} is based on the study of
properties of the Rauzy graph and its behaviour with respect to
the operation of mirror image on the language of the infinite
word. It turns out that the Rauzy graphs of words reaching the
upper bound in our estimate of palindromic complexity must have a
very special form.

%%%%%%%%%%%%%%%%%%%%%%%%%%%%%%%%%%%%%%%%%%%%%%%%%%%%%%%%%%%%%%%%%%%%%%%
\section{Acknowledgment}

The authors acknowledge the financial support of Czech Science
Foundation GA\v{C}R 201/05/0169 and the Ministry of Education of
the Czech Republic LC00602.

%%%%%%%%%%%%%%%%%%%%%%%%%%%%%%%%%%%%%%%%%%%%%%%%%%%%%%%%%%%%%%%%%%%%%%%

\end{document}